\documentclass[leqno]{article}

\newtheorem{theorem}{Theorem}
\newtheorem{lemma}[theorem]{Lemma}
\newtheorem{proposition}[theorem]{Proposition}
\newtheorem{definition}[theorem]{Definition}
\newtheorem{corollary}[theorem]{Corollary}

\newcommand{\begintheorem}{\addtocounter{equation}{1}\begin{theorem}}
\newcommand{\beginlemma}{\addtocounter{equation}{1}\begin{lemma}}
\newcommand{\beginproposition}{\addtocounter{equation}{1}\begin{proposition}}
\newcommand{\begindefinition}{\addtocounter{equation}{1}\begin{definition}}
\newcommand{\begincorollary}{\addtocounter{equation}{1}\begin{corollary}}

\begin{document}

\title{Some topics in complex and harmonic analysis, 3}

\author{Stephen William Semmes	\\
	Rice University		\\
	Houston, Texas}

\date{}

\maketitle

	Fix a positive integer $n$, and let $\mathcal{S}({\bf R}^n)$
denote the Schwartz class of complex-valued smooth functions on ${\bf
R}^n$ which are rapidly decreasing and whose derivatives are rapidly
decreasing.  For each pair of multi-indices $\alpha$, $\beta$, define
a seminorm $\|\cdot \|_{\alpha, \beta}$ on $\mathcal{S}({\bf R}^n)$ by
\begin{equation}
	\|\phi \|_{\alpha, \beta} =
  \sup \biggl\{\biggl|x^\alpha \, \frac{\partial^{d(\beta)}}{\partial x^\beta}
			\, \phi(x) \biggr| : x \in {\bf R}^n \biggr\},
\end{equation}
where $d(\beta)$ denotes the degree of the multi-index $\beta =
(\beta_1, \ldots, \beta_n)$, which is equal to $\beta_1 + \cdots +
\beta_n$.  One can characterize the Schwartz class $\mathcal{S}({\bf
R}^n)$ as the set of smooth functions $\phi$ on ${\bf R}^n$ such that
the seminorms $\|\phi\|_{\alpha, \beta}$ of $\phi$ are finite for all
multi-indices $\alpha$, $\beta$.

	This family of seminorms determines a topology on the Schwartz
class $\mathcal{S}({\bf R}^n)$ in a natural way, which can be
described by saying that the balls in $\mathcal{S}({\bf R}^n)$ with
respect to the seminorms $\|\cdot \|_{\alpha, \beta}$ form a sub-basis
for the topology.  A sequence $\{\phi_j\}_{j=1}^\infty$ of functions
in $\mathcal{S}({\bf R}^n)$ converges to another function $\phi \in
\mathcal{S}({\bf R}^n)$ if and only if $\|\phi_j - \phi\|_{\alpha,
\beta} \to 0$ as $j \to \infty$ for all multi-indices $\alpha$,
$\beta$.  A subset $E$ of $\mathcal{S}({\bf R}^n)$ is closed if and
only if for every sequence $\{\phi_j\}_{j=1}^\infty$ of elements of
$E$ which converges to some $\phi \in \mathcal{S}({\bf R}^n)$ in the
sense just defined we have that $\phi \in E$.

	For each nonnegative integer $l$, let $V_l$ be the linear
subspace of $\mathcal{S}({\bf R}^n)$ consisting of functions which
vanish up to order $l$.  When $l = 0$ we have that $\phi \in
\mathcal{S}({\bf R}^n)$ lies in $V_0$ exactly when $\phi(0) = 0$, and
in general we have that $\phi \in V_l$ when $\partial^{d(\beta)} /
\partial x^\beta \phi(0) = 0$ for all multi-indices $\beta$ such that
$d(\beta) \le l$.  Let us define $V_\infty$ to be the linear subspace
of $\mathcal{S}({\bf R}^n)$ consisting of functions which vanish to
infinite order at $0$, which is to say that the function and all of
its derivatives at $0$ are equal to $0$.  Equivalently, $V_\infty$
is equal to the intersection of all of the $V_l$'s.

	Let $\mathcal{P}({\bf R}^n)$ denote the vector space of
polynomials on ${\bf R}^n$, and for each nonnegative integer $l$ let
$\mathcal{P}_l({\bf R}^n)$ denote the linear subspace of
$\mathcal{P}({\bf R}^n)$ consisting of polynomials of degree less than
or equal to $l$.  Thus $\mathcal{P}({\bf R}^n)$ consists of all finite
linear combinations of monomials $x^\alpha = x_1^{\alpha_1} \cdots
x_n^{\alpha_n}$, and $\mathcal{P}_l({\bf R}^n)$ is the subspace
spanned by the monomials $x^\alpha$ where the multi-index $\alpha$ has
degree less than or equal to $l$.  For each positive integer $l$
there is a natural linear mapping from $\mathcal{S}({\bf R}^n)$
onto $\mathcal{P}_l({\bf R}^n)$ with kernel equal to $V_l$, which
sends a function $\phi \in \mathcal{S}({\bf R}^n)$ to its associated
Taylor polynomial of degree $l$ at $0$.  In particular $V_l$
has finite codimension in $\mathcal{S}({\bf R}^n)$ for all $l$.

	It is easy to see that $V_l$ is a closed linear subspace of
$\mathcal{S}({\bf R}^n)$ for all nonnegative integers $l$, and that
$V_\infty$ is a closed subspace of $\mathcal{S}({\bf R}^n)$ too.  Let
$W$ denote the linear subspace of $\mathcal{S}({\bf R}^n)$ consisting
of functions $\phi \in \mathcal{S}({\bf R}^n)$ which vanish on a
neighborhood of $0$, i.e., so that there is an open subset $U$ of
${\bf R}^n$ such that $0 \in U$ and $\phi(x) = 0$ for all $x \in U$.
Clearly $W \subseteq V_\infty$.  One can show that $V_\infty$ is
actually the closure of $W$.

	Let $\mathcal{S}'({\bf R}^n)$ denote the vector space of
tempered distributions on ${\bf R}^n$, which is to say the vector
space of continuous linear functionals on $\mathcal{S}({\bf R}^n)$.
More concretely, a tempered distribution $\lambda$ is a linear mapping
from $\mathcal{S}({\bf R}^n)$ to the complex numbers such that
\begin{equation}
	|\lambda(\phi)| \le C \, \sum_{j = 1}^p \|\phi\|_{\alpha_j, \beta_j}
\end{equation}
for some nonnegative real number $C$ and multi-indices $\alpha_1,
\beta_2, \ldots, \alpha_m, \beta_p$.  As usual we can define
derivatives of a tempered distribution $\lambda$ by applying $\lambda$
to the corresponding derivatives of a function $\phi \in
\mathcal{S}({\bf R}^n)$, and where we multiply the result by $-1$ if
the order of differentiation is odd.  This reduces to ordinary
differentiation when $\lambda(\phi)$ is defined by integrating $\phi$
times a smooth function of polynomial growth whose derivatives also
have polynomial growth, by integration by parts.

	Let us say that a tempered distribution $\lambda$ is supported
at the origin if $\lambda(\phi) = 0$ for all $\phi \in
\mathcal{S}({\bf R}^n)$ which vanish on a neigborhood of the origin.
This is equivalent to saying that $\lambda(\phi) = 0$ for all $\phi
\in V_\infty$.  One can check that if $\lambda$ is a tempered
distribution which is supported at the origin, then there is a
nonnegative integer $l$ such that $\lambda(\phi) = 0$ for all $\phi
\in V_l$.

	A basic example of such a tempered distribution is the Dirac
delta function at $0$, which sends $\phi \in \mathcal{S}({\bf R}^n)$
to $\phi(0)$.  Similarly derivatives of the Dirac delta function at
$0$ are tempered distributions which are supported at $0$.  These can
be defined as the tempered distributions which send a function $\phi
\in \mathcal{S}({\bf R}^n)$ to a specific derivative of $\phi$ at $0$.
Conversely, one can show that every tempered distribution on ${\bf R}^n$
supported at the origin is a finite linear combination of derivatives
of the Dirac delta function at $0$.

	If $\lambda$ is a tempered distribution on ${\bf R}^n$, then
the Fourier transform $\widehat{\lambda}$ of $\lambda$ is the tempered
distribution defined by
\begin{equation}
	\widehat{\lambda}(\phi) = \lambda(\widehat{\phi}),
\end{equation}
where $\widehat{\phi}$ denotes the classical Fourier transform
of $\phi \in \mathcal{S}({\bf R}^n)$, as an integrable function.
This makes sense as a tempered distribution because the
Fourier transform defines a continuous linear mapping from
$\mathcal{S}({\bf R}^n)$ into itself.  The Fourier transform
of tempered distributions has many of the same properties of the
classical Fourier transform on more regular functions, as one can
verify by reducing to the Fourier transform applied to functions
in the Schwartz class.

	Let $h(x)$ be a continuous function on ${\bf R}^n$ which has
polynomial growth in the sense that $|h(x)| \le C \, (1 + |x|^m)$ for
some nonnegative real number $C$, nonnegative integer $m$, and all $x
\in {\bf R}^n$.  We can use $h$ to define a tempered distribution
$\lambda$ by putting $\lambda(\phi)$ equal to the integral of $\phi$
times $h$ on ${\bf R}^n$ for all $\phi \in \mathcal{S}({\bf R}^n)$.
Let us assume that $h$ is harmonic, which one can define in the sense
of tempered distributions by saying that $\lambda (\Delta \phi) = 0$
for all $\phi \in \mathcal{S}({\bf R}^n)$, where $\Delta$ is the usual
Laplacian,
\begin{equation}
	\Delta = \sum_{j=1}^n \frac{\partial^2}{\partial x_j^2}.
\end{equation}

	In terms of the Fourier transform we have that
$\widehat{\lambda}(\psi_1) = 0$ whenever $\psi_1(\xi) = |\xi|^2 \,
\psi(xi)$ and $\psi \in \mathcal{S}({\bf R}^n)$.  In particular it
follows that $\widehat{\lambda}$ is supported at the origin.
As before, we may conclude that $\widehat{\lambda}$ is a linear
combination of derivatives of Dirac delta functions at the origin.
This implies that $h$ is a polynomial.

\end{document}